\documentclass[12pt,reqno]{amsart}

\usepackage{diagrams,amssymb, amsmath, amscd}

\diagramstyle[height=2em,width=2em]
\newarrow{Epi}----{>>}
\newarrow{Mono}>--->
\newarrow{Iso}>---{>>}
\newarrow{Mapsto}|--->
\newarrow{Igual}=====
\newarrow{Dashto}{}{dash}{}{dash}>

\newcommand{\Z}{{\mathbb Z}}

\newcommand{\C}{{\mathbb C}}
\newcommand{\F}{{\mathbb F}}

\newcommand{\tb}{\overline{T}}

                     \newcommand{\Hom}{{\rm Hom}}

\newcommand{\map}{{\rm map}}

\newcommand{\im}{{\rm Im}}

\newcommand{\A}{\ifmmode{\mathcal{A}}\else${\mathcal{A}}$\fi}
\newcommand{\K}{\ifmmode{\mathcal{K}}\else${\mathcal{K}}$\fi}
\newcommand{\U}{\ifmmode{\mathcal{U}}\else${\mathcal{U}}$\fi}
\newcommand{\T}{\ifmmode{\mathcal{T}}\else${\mathcal{T}}$\fi}

\newtheorem{theorem}{Theorem}[section]
\newtheorem{proposition}[theorem]{Proposition}
\newtheorem{corollary}[theorem]{Corollary}
\newtheorem{lemma}[theorem]{Lemma}

\newtheorem{remark}[theorem]{Remark}
\newtheorem{example}[theorem]{Example}

\title[Deconstructing Hopf spaces]{Deconstructing Hopf spaces}

\author{Nat\`{a}lia Castellana}
\author{Juan A. Crespo}
\author{J\'er\^{o}me Scherer}

\thanks{The authors are partially
supported by MEC grant MTM2004-06686 and the third author by the
program Ram\'on y Cajal, MEC, Spain.}

\begin{document}


\begin{abstract}
We characterize Hopf spaces with finitely generated cohomology as
an algebra over the Steenrod algebra. We ``deconstruct" the
original space into an $H$-space $Y$ with finite mod $p$
cohomology and a finite number of $p$-torsion Eilenberg-Mac Lane
spaces. We give a precise description of homotopy commutative
$H$-spaces in this setting.
\end{abstract}


\maketitle



\section*{Introduction}
\label{sec intro}
Since their introduction in the 50's by Serre, $H$-spaces have
produced some of the most beautiful results in Algebraic Topology.
Some examples are Adams' solution of the Hopf invariant one
conjecture~\cite{MR0141119}, the criminal of Hilton-Roitberg
\cite{MR0246320}, the construction of $DI(4)$ by Dwyer and
Wilkerson \cite{MR1161306}, the recent proof that a finite loop
space is of the homotopy type of a manifold \cite{MR2079597}, and
the new example of a finite loop space in \cite{ABGP}.

The structure of finite $H$-spaces is rather well understood. In
one of the most important articles on finite $H$-spaces,
\cite{MR38:6592}, Hubbuck shows that there are no other finite
connected homotopy commutative $H$-spaces than products of
circles, which was proved for compact Lie groups by James. It was
not until the early $90$'s that this result was extended by Slack
to $H$-spaces with finitely generated mod $2$ cohomology. With the
aid of secondary operations, he shows in~\cite{MR92k:55015} that
such homotopy commutative $H$-spaces are products of circles and
other Eilenberg-Mac Lane spaces. In fact, using the modern
techniques of Lannes' $T$ functor, Broto et al. obtain
in~\cite{MR2002b:55015} a structure theorem for all $H$-spaces
with noetherian mod $p$ cohomology. They ``deconstruct" such an
$H$-space into mod $p$ finite ones and copies of $K(\Z/p^r, 1)$
and $\C P^\infty$ in a functorial way. Recall that an $H$-space is
said to be mod $p$ finite if it is $p$-complete with finite mod
$p$ cohomology, which we denote simply by $H^*(-)$.

Our goal is to extend such results to an even larger class of
infinite dimensional spaces and understand which are the basic
pieces permitting to reconstruct the original $H$-space.

Natural examples of $H$-spaces, arising in connection with those
which are finite, are their Postnikov sections and connected
covers. The mod $p$ cohomology of the $n$-connected cover of a
finite $H$-space is not finite in general, but is finitely
generated as an algebra over the Steenrod algebra $\A_p$ (we refer
to the article~\cite{CCS4} for more details). Up to
$p$-completion, the basic examples of $H$-spaces satisfying this
cohomological finiteness condition are all finite $H$-spaces, and
Eilenberg-Mac Lane spaces of type $K(\Z/p^r,n)$ and
$K(\Z_{p^\infty}, n)$.

We show that one can deconstruct any such $H$-space in terms of
these basic $H$-spaces. We call an $H$-space which has only
finitely many non-trivial homotopy groups an {\it $H$-Postnikov
piece}.
\medskip

{\bf Theorem~\ref{estructura-finitamente-generados-sobre-A}.}
\noindent {\it Let $X$ be a connected $H$-space such that $H^*(X)$
is  a finitely generated algebra over the Steenrod algebra. Then
$X$ is the total space of an $H$-fibration
$$
F \rTo X \rTo Y,
$$
where $Y$ is an $H$-space with finite mod $p$ cohomology and $F$
is a $p$-torsion $H$-Postnikov piece whose homotopy groups are
finite direct sums of copies of cyclic groups $\Z/p^r$ and
Pr\"ufer groups $\Z_{p^\infty}$.}

\medskip

The above fibration behaves well with respect to loop structures
and a similar result holds for loop spaces. Our deconstruction
theorem enables us to reduce questions on infinite dimensional
$H$-spaces to finite ones. For instance, we use this technique to
give a generalization of Hubbuck's Torus Theorem.
\medskip

{\bf Corollary~\ref{Corolario1}.}
\noindent {\it Let $X$ be a connected homotopy commutative
$H$-space such that the mod $2$ cohomology $H^*(X)$ is finitely
generated as an algebra over the Steenrod algebra $\A_2$. Then, up
to $2$-completion, $X$ is homotopy equivalent to $(S^1)^n\times
F$, where $F$ is a connected $2$-torsion $H$-Postnikov piece.}

\medskip

The above splitting is not a splitting of $H$-spaces in general,
think of $S^1 \times K(\Z/2, 2)$, \cite[Section~1.4]{MR55:13416}.
When $H^*(X)$ is finitely generated as an algebra, we get back
Slack's result~\cite{MR92k:55015}, as well as their generalization
by Lin and Williams in \cite{MR92b:55011}.

\medskip

The arguments to prove our main theorem are the following. When
$H^*(X)$ is finitely generated over $\A _p$, we show in
Lemma~\ref{fg in the filtration} that the unstable module of
indecomposable elements $QH^*(X)$ belongs to some stage $\U_n$ of
the Krull filtration in the category of unstable modules. This
filtration has been studied in~\cite{MR95d:55017} by Schwartz and
then used in~\cite{MR2002k:55043} in order to prove Kuhn's
non-realizability conjecture~\cite{MR96i:55027}.

The stage $\U_0$ of the Krull filtration is particularly
interesting since it consists exactly of all locally finite modules
(direct limits of finite modules). In fact, the condition that $QH^*(X)$ is
locally finite is  equivalent to requiring that the loop space
$\Omega X$ be $B\Z/p$-local, i.e. the space of pointed maps $\map
_*(B\Z/p, \Omega X)$ is contractible, see \cite[Prop
3.2]{MR92b:55004} and \cite[Proposition~6.4.5]{MR95d:55017}.

We extend this topological characterization to $H$-spaces $X$ with
$QH^*(X)\in \U _n$. We use the standard notation $T_V$ for Lannes'
$T$ functor and say that $H^*(X)$ is of finite type if $H^n (X)$
is a finite $\F_p$-vector space for any integer $n\geq 0$.

\medskip

{\bf Theorem~\ref{equivalence}.}
\noindent {\it Let $X$ be a connected $H$-space such that $T_V
H^*(X)$ is of finite type for any elementary abelian $p$-group
$V$. Then $QH^*(X)$ is in $\mathcal U_n$ if and only if
$\Omega^{n+1}X$ is $B\Z/p$-local.}

\medskip

We apply now Bousfield's results on the Postnikov-like tower
associated to the $B\Z/p$-nullification functor $P_{B\Z/p}$
(relying on his ``Key Lemma", \cite[Chapter~7]{B1}). They enable
us to reconstruct those $H$-spaces such that $\Omega^{n+1}X$ is
$B\Z/p$-local from $P_{B\Z/p} X$ in a finite number of principal
$H$-fibrations over $p$-torsion Eilenberg-Mac Lane spa\-ces. When
$n=0$, we recover the results of Broto et al. from
\cite{BC,MR2002g:55016,MR2002b:55015}.

\medskip

{\bf Theorem~\ref{Teorema de estructura1}.}
\noindent {\it Let $X$ be an $H$-space such that $T_V H^*(X)$ is
of finite type for any elementary abelian $p$-group $V$. Then
$QH^*(X)$ is in $\mathcal U_n$ if and only if $X$ is the total
space of an $H$-fibration
$$
F \rTo X \rTo P_{B\Z/p} X
$$
where $F$ is a $p$-torsion $H$-Postnikov piece whose homotopy
groups are finite direct sums of copies of cyclic groups $\Z/p^r$
and Pr\"ufer groups $\Z_{p^\infty}$ concentrated in degrees from
$1$ to $n+1$. }

\medskip

It is worthwhile to mention that working with $H$-spaces is
crucial as illustrated by the space $BS^3$. Its loop space is
$B\Z/p$-local, but the fiber of the nullification map has
infinitely many non-trivial homotopy groups (see Dwyer's
computations \cite[Theorem~1.7, Lemma~6.2]{MR97i:55028}).

\medskip


\noindent
{\bf Acknowledgement.} Most of this work has been done
in the coffee room of the Maths Department at the UAB. We would
like to thank Alfonso Pascual for his generosity. We warmly thank
Carles Broto for his questions which regularly opened new
perspectives and Jesper Grodal for many useful comments. Finally,
we are indebted to the referee for suggesting to use
Andr\'e-Quillen cohomology to fix an earlier proof of
Theorem~\ref{complement}.

\section{Lannes' $T$ functor}
\label{section Lannes T functor}

Lannes' $T$ functor, \cite{MR90h:55027}, was designed as a tool to
compute the cohomology of mapping spaces with source $BV$, the
classifying space of an elementary abelian $p$-group $V$. It was
used also by Lannes to give an alternative proof of Miller's
Theorem on the Sullivan's conjecture.

Let $\mathcal U$ (resp. $\mathcal K$) be the category of unstable
modules (resp. algebras) over the Steenrod algebra. The functor
$T_V$ is the left adjoint of  $-\otimes H^*(BV)$ in $\mathcal U$
and $\mathcal K$, where $V$ is an elementary abelian $p$-group.
The left adjoint of $-\otimes \widetilde{H}^*(BV)$ is called the
reduced $T$ functor and denoted by $\tb_V$. For each unstable
module $M\in { \mathcal U}$, we have a splitting of modules over
the Steenrod algebra $T_V M=M\oplus \tb_V M$. We will use $T$ to
denote $T_{\Z/p}$ and $\tb$ to denote $\tb_{\Z/p}$.

If $M=H^*(X)$, the evaluation map $BV\times \map (BV,X)\rightarrow
X$ induces by adjunction a map $\lambda_V:T_V H^*(X) \rTo H^*(\map
(BV,X))$ of unstable algebras over $\A_p$, which is often an
isomorphism. When working with $H$-spaces, it is often handy to
deal with the pointed mapping space instead of the full mapping
space.

\begin{proposition}
\label{T-para-H-spacios-calcula}
Let $X$ be an $H$-space such that $H^*(X)$ is of finite type.
Then, $T_V H^*(X)$ is of finite type if and only if $H^*(\map
_*(BV,X))$ is of finite type. In this case $ T_V H^*(X) \cong
H^*(\map (BV,X))$, as algebras over $\A_p$.
\end{proposition}

\begin{proof}
If $X$ is an $H$-space, then $\map (BV,X)$ is again an $H$-space,
and so is the connected component $\map  (BV,X)_c$ of the constant
map (see \cite{MR55:13416}). Moreover, when $X$ is connected, all
connected components of the mapping space have the same homotopy
type. Since the evaluation map is an $H$-map and has a section,
there is a splitting
$$
\map (BV,X) \simeq X \times \map _*(BV,X)\,.
$$
which allows to work with the pointed mapping space. By
\cite[Theorem~1.5]{Miller} there is a weak equivalence $\map
_*(BV,X) \simeq \map _*(BV,\hat{X}_p)$ for any elementary abelian
$p$-group $V$. Since $X$ is $p$-good, we can work with
$\hat{X}_p$. Now $T_V H^*(\hat{X}_p) \cong H^*(\map
(BV,\hat{X}_p))$ by \cite[Proposition 3.4.4]{MR93j:55019}.
\end{proof}

When $X$ is connected, the evaluation $\map (BV,X) \to X$ is a
homotopy equivalence if $T_V H^*(X)\cong H^*(X)$ (for finite
spaces, this is the Sullivan conjecture proved by Miller
\cite[Theorems A,C]{Miller}). Actually, spaces for which this
happens have been cohomologically characterized by Lannes and
Schwartz in \cite{MR90h:55027}: their mod $p$ cohomology is
locally finite.

When one restricts the evaluation map to the connected component
of the constant map, the module of indecomposable elements
$QH^*(X)$ comes into play as observed by Dwyer and Wilkerson in
\cite[Proposition 3.2]{MR92b:55004} (see also \cite[3.9.7 and
6.4.5]{MR95d:55017}).

\begin{proposition}
\label{cor-DW}
Let $X$ be a connected $H$-space of finite type. Then $QH^*(X)$ is
a locally finite $\A_p$-module if and only if $\map _*(BV,\Omega
X)$ is contractible for some elementary abelian $p$-group $V$.
\end{proposition}

\begin{proof}
Since $\hat X_p$ is a connected $p$-complete $H$-space, $QH^*(X)$
is a locally finite $\A_p$-module if and only if $\map _*(BV, \hat
X_p)$ is homotopically discrete for any elementary abelian
$p$-group $V$ by \cite[Proposition~3.2]{MR92b:55004} and
\cite[Proposition~6.4.5]{MR95d:55017}.

The weak equivalence $\map _*(BV,X) \simeq \map _*(BV,\hat X_p)$
given by \cite[Theorem~1.5]{Miller} shows that this is equivalent
to $\map _*(BV,\Omega X)$ being contractible, i.e. the loop space
$\Omega X$ is $BV$-local.
\end{proof}

\section{The Krull filtration of $\U$}
\label{Krull}

In \cite{MR2002k:55043}, Schwartz proves  the ``strong realization
conjecture" extending his previous results
from~\cite{MR99j:55019}. This conjecture, given by Kuhn in
\cite{MR96i:55027}, states that if the cohomology of a space lies
in some stage of the Krull filtration of the category $\U$ of
unstable modules, then it must be locally finite. The Krull
filtration $\U_0 \subset \U_1 \subset \dots$ is defined
inductively, see \cite[Section~6.2]{MR95d:55017}. It starts with
the full subcategory $\U_0$ of $\U$ of locally finite unstable
modules and the modules in $\U_n$ can be characterized as follows
by means of the functor~$\tb$:

\begin{theorem}\cite[Theorem~6.2.4]{MR95d:55017}
\label{TKrull}
Let $M$ be an unstable module. Then $M \in \U_n$ if and only if
$\tb ^{n+1}M=0$. \hfill \qed
\end{theorem}

The proof of Kuhn's conjecture by Schwartz shows that under the
usual finiteness conditions the cohomology of a space either lies
in $\U_0$ or it is not in any $\U _n$. Instead of looking at when
the full cohomology of a space is in $\U_n$, we will study the
module of the indecomposable elements $QH^*(X)$. The Krull
filtration induces a filtration of the category of $H$-spaces by
looking at those $H$-spaces $X$ for which $QH^*(X) \in \U_n$.
There exist many spaces lying in each layer of this filtration,
the most obvious ones being Eilenberg-Mac\,Lane spaces.

\begin{example}
\label{E-M spaces are in Tn}
The module $QH^*(K(\Z/2, n+1))$ is isomorphic to the suspension of
the free unstable module $F(n)$ on one generator in degree~$n$. In
particular, the formula $\tb F(n)=\oplus_{0\leq i \leq n-1} F(i)$
(see \cite[Lemma 3.3.1]{MR95d:55017}) yields that $QH^*(K(\Z/2,
n+1))\in \U_n$. For an odd prime $p$, an easy induction argument
based on the next lemma shows that $QH^*(K(\Z/p, n+1))\in \U_n$ as
well.

More generally, let $G$ be any abelian discrete group such that
$H^*(K(G,n+1))$ is of finite type. Then $QH^*(K(G,n+1))\in \U_{n}$
(see \cite[Section 9.8]{MR95d:55017} for the explicit computations
of the $T$ functor).
\end{example}

{}From the above example, we see that the filtration is not
exhaustive, since the infinite product $\prod_{n \geq 1} K(\Z/p,
n)$ does not belong to any stage of the filtration. As we explain
in Example~\ref{BU} nor does~$BU$.

Next lemma shows, by means of the reduced $T$ functor, how
$QH^*(X)$ is related to $QH^*(\map _*(B\Z/p,X) )$.

\begin{lemma}
\label{T-de-los-indescomponibles}
Let $X$ be an $H$-space such that $T H^*(X)$ is of finite type.
Then there is an isomorphism $\tb QH^*(X) \cong QH^*(\map
_*(B\Z/p,X) )$.
\end{lemma}

\begin{proof}
Under such assumptions, Proposition~\ref{T-para-H-spacios-calcula}
applies and we know that the $T$ functor computes the cohomology
of the mapping space. Thus $QT H^*(X)$ is isomorphic to
$$
QH^*(\map (B\Z/p, X)) \cong Q \bigr( H^*(\map _*(B\Z/p, X))
\otimes H^*(X) \bigl)
$$
Since $T$ commutes with taking indecomposable elements \cite[Lemma
6.4.2]{MR95d:55017}, it follows that $T Q H^*(X) \cong
QH^*(X)\oplus QH^*(\map _*(B\Z/p,X))$. This is equivalent to $\tb
QH^*(X)\cong QH^*(\map  _*(B\Z/p,X))$.
\end{proof}

We end the section with a proposition which will allow us to
perform an induction in the Krull filtration. Observe that Kuhn's
strategy to move in the Krull filtration is to consider the
cofiber of the inclusion $X \rightarrow \map (B\Z/p, X)$ in the component of the constant map, see
\cite{MR96i:55027}. In our context,
Lemma~\ref{T-de-los-indescomponibles} suggests to use the fiber of
the evaluation $\map (B\Z/p, X) \rightarrow X$.

\begin{proposition}
\label{H*-y-H*map}
Let $X$ be an $H$-space with $T H^*(X)$ of finite type. Then, for
$n\geq 1$, $QH^*(X)\in \U_ n$ if and only if $QH^*(\map _*(B\Z/p,
X))$ is in $\U_{n-1}$.
\end{proposition}

\begin{proof}
By Theorem~\ref{TKrull}, the unstable module $QH^*(X)$ belongs to
$\U_n $ if and only if $\tb ^{n+1}QH^*(X) =0 $. By
Lemma~\ref{T-de-los-indescomponibles}, $\tb^n \tb QH^*(X) \cong
\tb ^n(QH^*(\map _*(B\Z/p, X))$ , and we obtain that $QH^*(X) \in
\mathcal U _{n}$ if and only if $QH^*(\map _*(B\Z/p,X)) \in
\mathcal U _{n-1}$.
\end{proof}

By repeatedly applying the previous proposition, one can give a
more geometrical formulation to the condition $QH^*(X) \in \U_ n$.
This happens if and only if the pointed mapping space out of an
$(n+1)$-fold smash product $\map _*(B\Z/p \wedge \dots \wedge
B\Z/p, X)$ is homotopically discrete.
\section{Bousfield's $B\Z/p$-nullification filtration.}\label{local}
The plan of this section follows the preceding one step-by-step,
replacing the algebraic filtration defined with the module of
indecomposables by a topological one. Recall (cf. \cite{Dror})
that, given a pointed connected space $A$, a space $X$ is
$A$-\emph{local} if the evaluation at the base point in $A$
induces a weak equivalence of mapping spaces $\map(A, X) \simeq
X$. When $X$ is an $H$-space, it is sufficient to require that the
pointed mapping space $\map_*(A, X)$ be contractible.

Dror-Farjoun and Bousfield have constructed a localization functor
$P_{A}$ from spaces to spaces together with a natural
transformation $l: X \rightarrow P_{A}X$ which is an initial map
among those having an $A$-local space as target (see \cite{Dror}
and \cite{MR57:17648}). This functor is known as the
$A$-nullification. It preserves $H$-space structures since it
commutes with finite products. Moreover, when $X$ is an $H$-space,
the map $l$ is an $H$-map and its fiber is an $H$-space.

Bousfield has determined the structure of the fiber of the
nullification map $l:X \rightarrow P_{A}X$ under certain
assumptions on~$A$. We are interested in the situation in which
$A$ an iterated suspension of~$B\Z/p$.

\begin{theorem}\label{MT1}\cite[Theorem~7.2]{B2}
Let $n \geq 1$ and $X$ be a connected $H$-space such that
$\Omega^n X$ is $B\Z/p$-local. The homotopy fiber of the
localization map $X \rightarrow P_{\Sigma^{n-1}B\Z/p} X$ is then
an Eilenberg-Mac\,Lane space $K(P,n)$ where $P$ is an abelian
$p$-torsion group (possibly infinite). \hfill{\qed}
\end{theorem}

As mentioned by Bousfield in \cite[p. 848]{B2}, an inductive
argument allows to obtain a precise description of the fiber of
the $B\Z/p $-nullifica\-tion map for $H$-spaces for which some
iterated loop space is local.

\begin{theorem}\label{MT3}
Let $n\geq 0$ and $X$ be a connected $H$-space such that $\Omega^n
X$ is $B\Z/p$-local. Then there is an $H$-fibration
$$
F \rTo X \rTo P_{B\Z/p} X,
$$
where $F$ is a $p$-torsion $H$-Postnikov piece whose homotopy
groups are concentrated in degrees from $1$ to $n$. \hfill{\qed}
\end{theorem}

We introduce a ``nullification filtration" by looking at those
$H$-spaces $X$ such that the iterated loop space $\Omega^n X$ is
$B\Z/p$-local. The example of the Eilenberg-Mac Lane spaces shows
that there are many spaces living in each stage of this filtration
as well, compare with Example~\ref{E-M spaces are in Tn}.

\begin{example}
\label{E-M spaces are in Topn}
Let $G$ be an abelian discrete group with non-trivial mod $p$
cohomology. Then, the Eilenberg-Mac Lane space $K(G, n)$ enjoys the
property that its $n$ fold iterated loop space is $B\Z/p$-local
(it is even discrete). The infinite product $\prod_{n\geq 1}
K(\Z/p, n)$ does not live in any stage of this topological
filtration.
\end{example}

Another source of examples of spaces in this filtration is
provided by connected covers of finite $H$-spaces.

\begin{example}\label{recubridores y lazos}
Let $X$ be a finite connected $H$-space and consider its
$n$-connected cover $X\langle n \rangle$. Then
$\Omega^{n-1}(X\langle n \rangle)$ is $B\Z/p$-local.
\end{example}

For a connected $H$-space $X$ such that $\Omega^n X$ is
$B\Z/p$-local, the study of the homotopy type of  $\map
_*(B\Z/p,X)$ is drastically simplified by Theorem~\ref{MT3}, since
this space is equivalent to $\map _*(B\Z/p,F)$ where $F$ is a
Postnikov piece, as we explain in the proof below. A complete
study of the $B\Z/p$-homotopy theory of such $H$-spaces is
undertaken in \cite{CCS2}.

We prove now the topological analogue of
Lemma~\ref{T-de-los-indescomponibles}.

\begin{proposition}
\label{reduction}
Let $X$ be a  connected $H$-space such that $\Omega^n X$ is
$B\Z/p$-local, then $\Omega^{n-1} \map _*(B\Z/p,X)$ is
$B\Z/p$-local.
\end{proposition}

\begin{proof}
Under the hypothesis that $\Omega^{n}X$ is $B\Z/p$-local,
Theorem~\ref{MT3} tells us that we have a fibration
$$
F \rTo X \rTo P_{B\Z/p}X,
$$
where $F$ is a $p$-torsion Postnikov system with homotopy
concentrated in degrees from $1$ to $n$. Thus,  $\map  _*(B\Z/p,
X)\simeq \map  _*(B\Z/p, F)$ because $P_{B\Z/p}X$ is a
$B\Z/p$-local space. Now, $\Omega ^{n-1}\map  _*(B\Z/p, F)$ is
$B\Z/p$-local (in fact, it is a homotopically discrete space) and
thus so is $\Omega ^{n-1} \map  _*(B\Z/p, X)$.
\end{proof}

\section{Infinite loop spaces}
\label{sec infinite}
In order to compare the topological with the algebraic filtration,
one of the key ingredients comes from the theory of infinite loop
spaces. In this section we explain when a pointed mapping space
$\map _*(A, X)$ is an infinite loop space, but we are of course
specially interested in the case when $A$ is $B\Z/p$. We make use
of Segal's techniques of $\Gamma$-spaces and follow his notation
from \cite{MR50:5782}, which is better adapted to our needs than
that of Bousfield and Friedlander, see \cite{MR80e:55021}. Recall
that the category $\Gamma$ is the category of finite sets, and a
morphism $\theta: S \rightarrow T$ between two finite sets is a
partition of a subset of $T$ into $|S|$ disjoint subsets $\{ \,
\theta(\alpha) \, \}_{\alpha \in S}$. A $\Gamma$-space is a
contravariant functor from $\Gamma$ to the category of spaces with
some extra conditions. We first construct a covariant functor
$A_\bullet: \Gamma \rightarrow Spaces$ for any pointed space $A$
by setting $A_n = A^n$ (so in particular $A_0 = *$) and a morphism
$\theta: [n] \rightarrow [m]$ induces the map $\theta_*: A^n
\rightarrow A^m$ sending $(a_1, \dots, a_n)$ to the element $(b_1,
\dots, b_m)$ with $b_j = a_i$ if and only if $j \in \theta(i)$ and
$b_j = *$ otherwise.

Hence, we get a contravariant functor for any pointed space $X$ by
taking the pointed mapping space $\map _*(-, X)$. For
$\map_*(A_\bullet, X)$ to be a $\Gamma$-space one needs to check
it is {\it special}, i.e. the $n$ inclusions $i_k: [1] \rightarrow
[n]$ sending $1$ to $k$ must induce a weak equivalence $\map
_*(A^n, X) \rightarrow \map _*(A, X)^n$.

\begin{lemma}
\label{Gamma}
Let $A$ and $X$ be pointed spaces and assume that the inclusion
$A^n \vee A \hookrightarrow A^n \times A$ induces for any $n \geq
1$ a weak equivalence $\map _*(A^n \times A, X) \rightarrow \map
_*(A^n \vee A, X)$. Then, $\map _*(A_\bullet, X)$ is a
$\Gamma$-space. \hfill \qed
\end{lemma}

\begin{proposition}
\label{GammaH}
Let $A$ be a pointed connected space and $X$ an $H$-space. Assume
that  $\map _*(A, X)$ is $A$-local. Then $\map _*(A_\bullet, X)$
is a $\Gamma$-space.
\end{proposition}

\begin{proof}
The cofiber sequence $A^n \vee A \rightarrow A^n \times A
\rightarrow A^n \wedge A$ yields a fibration of pointed mapping
spaces
$$
\map _*(A^n \wedge A, X) \rTo \map _*(A^n \times A, X) \rTo \map
_*(A^n \vee A, X).
$$
By adjunction, the fiber $\map _*(A^n \wedge A, X) \simeq \map
_*(A^n, \map _*(A, X))$ is contractible since any $A$-local space
is also $A^n$-local ($A^n$ is $A$-cellular or use Dwyer's version
of Zabrodsky's Lemma in \cite[Proposition~3.4]{MR97i:55028}). Moreover, the
inclusion $A^n \vee A \rightarrow A^n \times A$ induces a
bijection on sets of homotopy classes $[A^n \times A, X]
\rightarrow [A^n \vee A, X]$ by \cite[Lemma~1.3.5]{MR55:13416}. Since
all components of these pointed mapping spaces have the same
homotopy type, we have a weak equivalence $\map _*(A^n \times A, X)
\simeq \map _*(A^n \vee A, X)$ and conclude by the preceding
proposition.
\end{proof}

\begin{theorem}
\label{infiniteloop}
Let $A$ be a pointed connected space and let $X$ be a loop space such
that $\map _*(A, X)$ is $A$-local. Then, $\map _*(A, X)$ is an
infinite loop space, and so is the corresponding connected
component $\map _*(A, X)_c$ of the constant map.
\end{theorem}

\begin{proof}
From the $\Gamma$-space structure constructed above we obtain
classifying spaces $B^n \map _*(A, X)$ and weak equivalences
$\Omega B^{n+1} \map _*(A, X) \simeq B^n \map _*(A, X)$ for any $n
\geq 1$. In our situation $X$ is a loop space, and so is the
mapping space $\map _*(A, X)$. Therefore, Segal's result
\cite[Proposition~1.4]{MR50:5782} applies and shows that $\map
_*(A, X)$ is equivalent to the loop space $\Omega B \map _*(A,
X)$.
\end{proof}

We specialize now to the case $A = B\Z/p$, where we can even say
more about the intriguing infinite loop space $\map _*(B\Z/p,
X)_c$. In the context of Proposition~\ref{superteoremazo} it will
turn out to be contractible.

\begin{proposition}
\label{vectorspace}
Let $X$ be a loop space such that $\map _*(B \Z/p, X)$ is
$B\Z/p$-local. Then all homotopy groups of the infinite loop space
$\map _*(B\Z/p, X)_c$ are $\Z/p$-vector spaces.
\end{proposition}

\begin{proof}
Since $\pi_n \map _*(B\Z/p, X)_c \cong [B\Z/p, \Omega^n X]$,
consider a map $B\Z/p \rightarrow \Omega^n X$. We claim that it is
homotopic to an $H$-map. Indeed, by \cite[Proposition~
1.5.1]{MR55:13416}, the obstruction lives in the set $[B\Z/p
\wedge B\Z/p, \Omega^n X]$, which is trivial since $\map _*(B
\Z/p, X)$ is $B\Z/p$-local. But any non-trivial $H$-map out of
$B\Z/p$ has order~$p$.
\end{proof}

\section{Structure theorems for $H$-spaces}
\label{section structure theorems}

The purpose of this section is to give an inductive description of
the $H$-spaces whose module of indecomposable elements lives in
some stage of the Krull filtration. This is achieved by comparing
this algebraic filtration with the topological one and by making
use of Bousfield's result \ref{MT3}.

\begin{proposition}
\label{teoremazo}
Let $X$ be an $H$-space such that $T_V H^*(X)$ is of finite type
for any ele\-mentary abelian $p$-group $V$. Assume that
$\Omega^nX$ is a $B\Z/p$-local space. Then $QH^*(X)~\in~\mathcal U
_{n-1}$.
\end{proposition}

\begin{proof}
We proceed by induction.  For $n=1$, assume that $\Omega X$ is
$B\Z/p$-local, that is, $\Omega \map _*(B\Z/p,X)_c\simeq *$. Then,
$\map _*(B\Z/p,X)$ is homotopically discrete since $\map
_*(B\Z/p,X)_c$  is so and all components of the mapping space have
the same homotopy type. Hence, $QH^*(\map _*(B\Z/p,X))=0$ and, by
Lemma  \ref{T-de-los-indescomponibles}, $QH^*(X)\in \mathcal
U_{0}$.

If $ n>1$, let $X$ be an $H$-space such that $\Omega^nX$ is
$B\Z/p$-local. We see by Proposition~\ref{reduction} that
$\Omega^{n-1}\map  _*(B\Z/p, X)_c$ is $B\Z/p$-local as well. Now,
$\map _*(B\Z/p,X)_c$ is an $H$-space such that $\Omega^{n-1}\map
_*(B\Z/p,X)_c$ is $B\Z/p$-local. Moreover, by
Proposition~\ref{T-para-H-spacios-calcula}, $T_V H^*(\map _*
(B\Z/p,X))$ is of finite type for any elementary abelian $p$-group
$V$. By induction hypothesis, $QH^*(\map  _*(B\Z/p, X)_c) \in
\mathcal U_{n-2}$. Since all components have the same homotopy
type, $QH^*(\map _*(B\Z/p, X)) \in \mathcal U_{n-2}$, and we
conclude that $QH^*(X) \in \mathcal U _{n-1}$  by
Proposition~\ref{H*-y-H*map}.
\end{proof}

\begin{proposition}
\label{superteoremazo}
Let $X$ be a connected $H$-space such that $T_V H^*(X)$ is of
finite type for any elementary abelian $p$-group $V$. Suppose that
$QH^*(X) \in \mathcal U_n$. Then $\Omega^{n+1}X$ is $B\Z/p$-local.
\end{proposition}

\begin{proof}
As in the proof of Proposition~\ref{T-para-H-spacios-calcula} we
can assume without loss of generality that $X$ is $p$-complete.

Let us proceed by induction. The case $n=0$ is given by
Proposition~\ref{cor-DW}. Now assume that the result is true for
$n-1$, and consider a space $X$ such that $QH^*(X)\in \U _n$.
Then, by Proposition~\ref{H*-y-H*map}, $QH^*(\map _*(B\Z/p, X))
\in \U_{n-1}$ and the induction hypothesis ensures that $\Omega^n
\map _*(B\Z/p,X)_c \simeq \map _*(B\Z/p, \Omega ^nX)$ is
$B\Z/p$-local. Apply now Theorem~\ref{infiniteloop} to deduce that
the space $\map _*(B\Z/p, \Omega ^nX)_c$ is an infinite loop
space, with a $p$-torsion fundamental group by
Proposition~\ref{vectorspace}.

These are precisely the conditions of McGibbon's main theorem in
\cite{MR97g:55013}: the $B\Z/p$-nullification of connected
infinite loop spaces with $p$-torsion fundamental group is trivial,
up to $p$-completion. Moreover, our infinite loop space is
$B\Z/p$-local, so
$$
(\map _*(B\Z/p, \Omega ^nX)_c)^{\wedge}_p \simeq \bigl(P_{B\Z/p}
(\map _*(B\Z/p, \Omega ^nX)_c) \bigr)^{\wedge}_p\simeq *
$$
As we assume that $X$ is $p$-complete, so are the loop space
$\Omega ^nX$ and the pointed mapping space $\map _*(B\Z/p, \Omega
^nX)_c$. Thus, we see that $\map _*(B\Z/p, \Omega ^nX)_c$ must be
contractible. Since all components of the pointed mapping space
have the same homotopy type as the component of the constant map,
we infer that $\map _*(B\Z/p, \Omega ^nX)$ is homotopically
discrete. Looping once again, one obtains finally a weak
equivalence $\map _*(B\Z/p, \Omega ^{n+1}X)\simeq *$, i.e. $\Omega
^{n+1}X$ is $B\Z/p$-local, as we wanted to prove.
\end{proof}

Let us sum up these two results in one single statement, which
extends widely Dwyer and Wilkerson's
\cite[Proposition~3.2]{MR92b:55004} when $X$ is assumed to be an
$H$-space.

\begin{theorem}
\label{equivalence}
Let $X$ be a connected  $H$-space such that $T_V H^*(X)$ is of finite type
for any elementary abelian $p$-group $V$. Then, $QH^*(X)$ is in
$\mathcal U_n$ if and only if $\Omega^{n+1}X$ is $B\Z/p$-local.
\hfill{\qed}
\end{theorem}

Combining these results with Theorem~\ref{MT3} (about the
nullification functor $P_{B\Z/p}$) enables us to give a
topological description of the $H$-spaces $X$ for which the
indecomposables $QH^*(X)$ live in some stage of the Krull
filtration.
Recall that the {\it Pr\"ufer group} $\Z_{p^\infty}$ is defined as
the union of all $\Z/p^n$, $n \geq 1$.
It is a $p$-torsion divisible abelian group.

\begin{theorem}\label{Teorema1}
Let $X$ be a connected $H$-space such that $T_V H^*(X)$ is of
finite type for any elementary abelian $p$-group $V$. Then
$QH^*(X)$ is in $\U_n$ if and only if $X$ fits into a principal
$H$-fibration
$$
K(P,n+1)\rTo X \rTo Y,
$$
where $Y$ is a connected $H$-space such that $QH^*(Y)\in \U_{n-1}$,
and $P$ is a $p$-torsion abelian group which is a finite direct
sum of copies of cyclic groups $\Z/p^r$ and Pr\"ufer groups
$\Z_{p^\infty}$.
\end{theorem}

\begin{proof}
Assume that $QH^*(X)$ is in $\U_n$. Let $F$ be the homotopy fiber
of the nullification map $X\rightarrow P_{\Sigma^{n}B\Z/p}(X)$. By
Theorem \ref{MT1}, $F\simeq K(P,n+1)$ where $P$ is an abelian
$p$-group. Moreover, the equivalence $\map _*(\Sigma^{n}B\Z/p,
K(P, n+1) ) \simeq \map _*(\Sigma^{n}B\Z/p, X)$ shows that the set
$$
\pi_{n} \map _*(B\Z/p, X) \cong \pi_0 \map _*(\Sigma^{n}B\Z/p, X)
\cong \Hom(\Z/p,P)
$$
is finite since all homotopy groups of $\map _*(B\Z/p, X)$ are
$p$-torsion and its cohomology is of finite type. Thus, $P$ is
isomorphic to a finite direct sum of copies of cyclic groups
$\Z/p^r$ and Pr\"ufer groups $\Z_{p^\infty}$ by Lemma \ref{finite
type p-groups}, which we prove at the end of the section.

We conclude by taking $Y=P_{\Sigma^{n}B\Z/p}(X)$. The cohomology
$H^*(Y)$ is of finite type since $H^*(K(P,n+1))$ and $H^*(X)$ are
of finite type, and so is $H^* (\map _*(BV, Y))$. Moreover, since
$\Omega^{n}Y \simeq P_{B\Z/p}(\Omega^{n}X)$ is $B\Z/p$-local,
Theorem~\ref{equivalence} implies that $QH^*(Y)\in \U_{n-1}$.
\end{proof}

Equivalently, one can reformulate this result by describing the
fiber of the $B\Z/p$-nullification map.

\begin{theorem}
\label{Teorema de estructura1}
Let $X$ be an $H$-space such that $T_V H^*(X)$ is of finite type
for any elementary abelian $p$-group $V$. Then, $QH^*(X)$ is in
$\mathcal U_n$ if and only if $X$ is the total space of an
$H$-fibration
$$
F \rTo X \rTo P_{B\Z/p} X,
$$
where $F$ is a $p$-torsion $H$-Postnikov piece whose homotopy
groups are finite direct sums of copies of cyclic groups $\Z/p^r$
and Pr\"ufer groups $\Z_{p^\infty}$ concentrated in degrees $1$ to
$n+1$. \hfill{\qed}
\end{theorem}

In other words, the only $H$-spaces such that $QH^*(X)$ lies in
$\U_n$ for some $n$ are the $B\Z/p$-local $H$-spaces, the
$p$-torsion Eilenberg-MacLane spaces introduced in
Example~\ref{E-M spaces are in Tn}, and extensions of the previous
type.

Recall that the $B\Z /p$-nullification of a loop space is again a
loop space. Moreover, by \cite[Lemma~3.A.3]{Dror}, the
nullification map is a loop map, and hence its homotopy fiber is
also a loop space. Thus we obtain automatically the following
result about loop spaces.

\begin{corollary}
\label{Teorema de estructura-loop spaces}
Let $X$ be a loop space such that $T_V H^*(X)$ is of finite type
for any elementary abelian $p$-group $V$. Then $QH^*(X)$ is in
$\mathcal U_n$ if and only if $X$ is the total space of loop
fibration
$$
F \rTo X \rTo P_{B\Z/p} X,
$$
where the loop space $F$ is a $p$-torsion Postnikov piece whose
homotopy groups are finite direct sums of copies of cyclic groups
$\Z/p^r$ and Pr\"ufer groups $\Z_{p^\infty}$ concentrated in
degrees from $1$ to $n+1$. \hfill{\qed}
\end{corollary}

If we restrict our attention to the case $n=0$ in
Theorem~\ref{Teorema de estructura1}, our result reproves in a
more conceptual way the theorems about $H$-spaces with locally
finite module of indecomposable elements given by Broto, Saumell
and the second named author
in~\cite{BC,MR2002g:55016,MR2002b:55015}.

What do we learn from our study about $H$-spaces which do not
belong to any stage of the filtration we have introduced in this
paper? From a cohomological point of view, such $H$-spaces have a
very large module of indecomposables. Let us discuss the example
of $BU$.

\begin{example}
\label{BU}
The module $QH^*(BU)$ is isomorphic to $\Sigma^2H^*(BS^1)$. In
particular, it is not a finitely generated  $\A_p$-module. A
computation of the value of the $T$ functor on this module can be
done using \cite[Section 9.8]{MR95d:55017} and shows that
$QH^*(BU)$ does not belong to any $\U_n$.

Therefore the Krull filtration for the indecomposables detects in
$BU$ the fact that the $B\Z/p$-nullification Postnikov-like tower
does not permit to deconstruct it into elementary pieces. In fact
$BU$ is $K(\Z/p, 2)$-local by a result of Mislin (see
\cite[Theorem~2.2]{MR57:7585}).
\end{example}

Finally, we prove  the  lemma about abelian $p$-torsion groups
which was used in the proof of Theorem~\ref{Teorema1}.


\begin{lemma}
\label{finite type p-groups}
Let $P$ be an abelian $p$-torsion group. If $\Hom(\Z/p,P)$ is
finite then $P$ is a finite direct sum of copies of cyclic groups
$\Z/p^r$ and Pr\"ufer groups $\Z_{p^\infty}$.
\end{lemma}

\begin{proof}
By Kulikov's theorem (see \cite[Theorem 10.36]{rotman}), $P$
admits a basic subgroup, which is a direct sum of cyclic groups.
It must be of bounded order since $\Hom(\Z/p,P)$ is finite, and a
result of Pr\"ufer (see \cite[Corollary 10.41]{rotman}) shows now
that this subgroup is a direct summand. Since the quotient is
divisible and $\Hom(\Z/p,P)$ is finite, $P$ is a finite direct sum
of copies of cyclic groups $\Z/p^r$ and Pr\"ufer groups
$\Z_{p^\infty}$.
\end{proof}

\section{Fibrations over Eilenberg-MacLane spaces}
\label{section EM}
In the next section we concentrate on $H$-spaces whose mod $p$
cohomology is finitely generated as an algebra over the Steenrod
algebra. Therefore we will need to establish a closure property
under certain $H$-fibrations.


\begin{theorem}
\label{complement}
Let $A$ be a finite direct sum of copies of cyclic groups $\Z/p^r$
and Pr\"ufer groups $\Z_{p^\infty}$, and $n \geq 2$. Consider an
$H$-fibration $F \xrightarrow{i} E \xrightarrow{\pi} K(A, n)$. If
$H^*(F)$ is a finitely generated $\A_p$-algebra, then so is
$H^*(E)$.
\end{theorem}

The proof relies mainly on the Eilenberg-Moore spectral sequence
for an $H$-fibration over an Eilenberg-Mac Lane space, a situation
studied by Smith in \cite[Chapter II]{MR0275435}. Following the
terminology used in \cite[Section~6]{MR0275435}, we say that a
sequence of (Hopf) algebras
$$
\F_p \rTo A \rTo B \rTo C \rTo \F_p
$$
is coexact if the morphism $A \rightarrow B$ is a monomorphism and
its cokernel $B//A$ is isomorphic to $C$ as a (Hopf) algebra.

\begin{proposition}
\label{coexact}
Let $n \geq 2$ and consider a non-trivial $H$-fibration $F
\xrightarrow{i} E \xrightarrow{\pi} K(A, n)$ where $A$ is either
$\Z/p$ or a Pr\"ufer group $\Z_{p^\infty}$. Then there is a
coexact sequence of algebras
$$
\F_p \rTo H^*(E)//\pi^* \rTo^{i^*} H^*(F) \rTo \Lambda \otimes S
\rTo \F_p\, ,
$$
where $\Lambda$ is an exterior algebra and $S \subseteq H^*(K(A,
n-1))$ is a Hopf subalgebra.
\end{proposition}

\begin{proof}
The $E_2$-term of the Eilenberg-Moore spectral sequence is given
by $E_2=\hbox{\rm Tor}_{H^*(K(A,n))}(H^*(E),\F_p)$ and converges
to $H^*(F)$. Since we deal with an $H$-fibration, an argument
based on the change of rings spectral sequence,
\cite[XVI.6.1]{MR1731415}, allows Smith to identify $E_2$ with
$H^*(E)//\pi^*\otimes \hbox{\rm
Tor}_{H^*(K(A,n))\backslash\backslash \pi^*}(\F_p,\F_p)$ as
algebras, see \cite[Theorem~2.4]{MR0275435}. Here
$H^*(K(A,n))\backslash\backslash \pi^*$ is the (Hopf algebra)
kernel of $\pi^*$, an unstable Hopf subalgebra of the abelian Hopf
algebra $H^*(K(A,n))$, as explained in
\cite[Remark~3.2]{MR0275435}.

By \cite[Proposition 7.3$^*$]{MR0275435}, $i^*: H^*(E)//\pi^*
\hookrightarrow H^*(F)$ is a monomorphism and the cokernel $R =
H^*(F)//i^*$ is described by a coexact sequence of Hopf algebras
$$
\F_p\rTo \Lambda \rTo  R \rTo
S \rTo \F_p,
$$
with $\Lambda$ an exterior algebra and $S \subseteq
H^*(K(A,n-1))$. Under our assumptions for $A$, $H^*(K(A,n-1))$ is
a free commutative algebra. Therefore the previous coexact
sequence splits, i.e. $R$ is isomorphic to $\Lambda \otimes S$ as
an algebra.
\end{proof}

To prove Theorem~\ref{complement}, we will show that the module of
indecomposable elements $QH^*(E)$ is a finitely generated unstable
module. When $p$ is an odd prime, the coexact sequence in the
proposition splits, which identifies $Q(H^*(E)//\pi^*)$ as a
submodule of $QH^*(F)$.

When $p=2$, the coexact sequence does not split. The functor $Q$
is not left exact and one is then naturally led to studying the
left derived functors of $Q$, i.e. Andr\'e-Quillen homology. Good
references are \cite[Chapter~7]{MR95d:55017}, \cite{MR1089001},
and \cite{Miller}. We write $H^Q_i(A)$ for the $i$-th derived
functor. This is an unstable module when $A$ is an unstable
algebra.

\begin{lemma}
\label{AQ}
Let $\Lambda$ be an exterior algebra over $\F_2$, which is
finitely generated as an algebra over the Steenrod algebra. Then
$H_1^Q(\Lambda)$ is a finitely generated unstable module.
\end{lemma}

\begin{proof}
In \cite[Section~10]{MR1089001}, Goerss identifies the first
Andr\'e-Quillen homology group $H_1^Q(\Lambda)$ with the
indecomposable elements of degree~$2$ in $\hbox{\rm
Tor}_\Lambda(\F_2,\F_2)$. As an $\F_2$-vector space it is
generated by the elements $[x|x]$ where $x$ runs through all
exterior generators of $\Lambda$. Since the Steenrod algebra acts
via the Cartan formula, it follows that $H_1^Q(\Lambda)$ is also a
finitely generated $\A_2$-module. Alternatively one could perform
this computation using the simplicial resolution given by the
symmetric algebra comonad described in \cite{MR95d:55017}.

\end{proof}

\noindent {\it Proof of Theorem~\ref{complement}.} Given a group
extension $A'\rightarrow A\rightarrow A$, there is a pullback
diagram of fibrations:
$$\begin{diagram}
E'& \rTo & E & \rTo & K(A'',n)
\\
\dTo & & \dTo_\pi & & \dIgual
\\
K(A',n) & \rTo & K(A,n) & \rTo & K(A'',n)
\\
\end{diagram}$$
If the statement is true for the first vertical fibration and the
top horizontal fibration, then it follows for $\pi$. Therefore, we
can assume that $A=\Z/p$ or $\Z_{p^\infty}$. If $\pi$ is
null-homotopic, the statement is obvious, and we can hence work
under the hypothesis of Proposition~\ref{coexact}.

Since $H^*(K(A,n))$ is finitely generated as algebra over
$\mathcal A_p$, so is its image $\im(\pi^*)\subseteq H^*(E)$.
Hence, to prove the theorem, it is enough to show that
$H^*(E)//\pi^*$ is a finitely generated $\A_p$-algebra, or
equivalently that the module of indecomposable elements
$Q(H^*(E)//\pi^*)$ is a finitely generated $\A_p$-module.

When $p$ is odd, the coexact sequence in Proposition~\ref{coexact}
splits (as algebras) because an exterior algebra is free
commutative. Hence $H^*(F) \cong H^*(E)//\pi^*\otimes \Lambda
\otimes S$ (compare with \cite[Theorem~5.7]{MR34:8406}). In
particular, $Q(H^*(E)//\pi^*) \subseteq QH^*(F)$. Since $\mathcal
U$ is a locally noetherian category,
\cite[Theorem~1.8.1]{MR95d:55017}, and $QH^*(F)$ is a finitely
generated $\A_p$-module, so is $Q(H^*(E)//\pi^*)$.

The case $p=2$ is less straightforward since an exterior algebra
is not free commutative, so that the coexact sequence in
Proposition~\ref{coexact} does not split in general. The inclusion
$A \subset B$ of a sub-Hopf algebra is not necessarily a
cofibration (seen as a constant simplicial object). However, when
$B$ is of finite type, it is always a free $A$-module by the
Milnor-Moore result \cite[Theorem~4.4]{MR0174052}). Therefore the
argument in \cite[Section~10]{MR1089001} applies and the homotopy
cofiber of the inclusion is weakly equivalent to $B//A$. Since
cofibrations of simplicial algebras induce long exact sequences in
Andr\'e-Quillen cohomology, we have in our situation an exact
sequence
$$
H_1^Q(\Lambda \otimes S) \rightarrow Q(H^*(E)//\pi^*) \rightarrow
QH^*(F).
$$
As $S$ is a free commutative algebra, $H_1^Q(\Lambda \otimes S)
\cong H_1^Q(\Lambda)$ by \cite[Lemma~4.10]{MR1089001}. The
previous lemma tells us that this is a finitely generated unstable
module. So is $QH^*(F)$, and we conclude since $\mathcal U$ is
locally noetherian. \hfill \qed

\begin{remark}
\label{pi1}
{\rm Theorem~\ref{complement} is actually true for fibrations over
an Eilen\-berg-Mac Lane space $K(A, 1)$ as well. The
Eilenberg-Moore spectral sequence converges in this case by work
of Dwyer, \cite{MR0394663}, and the above proof only needs minor
modifications.}
\end{remark}

\section{$H$-spaces with finitely generated cohomology over ${\mathcal A}_p$}
\label{section finitely generated}

We will assume in this section that $H^*(X)$ is finitely generated
as an algebra over the Steenrod algebra. Then, the
$B\Z/p$-nullification of $X$ is a mod $p$ finite $H$-space up to
$p$-completion, as we prove in Theorem~\ref{Teorema2}.

We show first that, under this finiteness condition, the
$H$-spaces considered in this section satisfy the hypothesis of
Theorem~\ref{equivalence} (they belong to some stage of the
filtration we study in this paper).

\begin{lemma}
\label{fg in the filtration}
Let $K$ be a finitely generated unstable ${\mathcal A}_p$-algebra.
Then there exists some integer $n$ such that the module of
indecomposables $QK$ belongs to $\U_n$. Moreover $T_V K$ is a
finitely generated unstable ${\mathcal A}_p$-algebra for any
elementary abelian group $V$.
\end{lemma}

\begin{proof}
First of all, $QK$ is a finitely generated module over ${\mathcal
A}_p$, i.e. it is a quotient of a finite direct sum of free
modules. Hence, there exists an epimorphism $\oplus_{i=1}^k
F(n_i)\twoheadrightarrow QK$. Since $\tb$ is an exact functor, it
follows that $\tb^m(QK)=0$, where $m$ is the largest of the
$n_i$'s, and so $QK \in \U_{m-1}$.

Moreover, $T_V$ commutes with taking indecomposables elements
\cite[Lemma 6.4.2]{MR2002k:55043}. Therefore, $Q(T_V K)$ is a
finitely generated unstable module. Then, the above discussion
shows that $T_V K$ is a finitely generated ${\mathcal
A}_p$-algebra.
\end{proof}

We can now state our main finiteness result. It enables us to
understand better the $B\Z/p$-nullification, which is the first
building block in our deconstruction process (Theorem~\ref{Teorema
de estructura1}).

\begin{theorem}\label{Teorema2}
Let $X$ be a connected $H$-space such that $H^*(X)$ is finitely
generated as algebra over the Steenrod algebra. Then, $P_{B\Z/p}
X$ is an $H$-space with finite mod $p$ cohomology.
\end{theorem}

\begin{proof}
By Lemma \ref{fg in the filtration}, there exists an integer $n$
such that $QH^*(X)$ lies in $\U_{n-1}$, so
Theorem~\ref{equivalence} applies and we know that $\Omega^nX$ is
$B\Z/p$-local.

We will show that if $H^*(X)$ is finitely generated as an algebra
over ${\mathcal A}_p$ and $\Omega^nX$ is $B\Z/p$-local, then
$H^*(P_{B\Z/p}X)$ is finitely generated as an algebra over ${\mathcal
A}_p$. We proceed by induction on $n$. When $n=0$ the statement is
clear. Assume the statement holds for $n-1$.

We know from Theorem~\ref{Teorema1} that there is a principal
$H$-fibration
$$
K(P,n) \rTo X \rTo P_{\Sigma^{n-1}B\Z/p}X,
$$
where $P$ is a $p$-torsion abelian group which is a finite direct
sum of copies of cyclic groups $\Z/p^r$ and Pr\"ufer groups
$\Z_{p^\infty}$. It follows from Theorem~\ref{complement} that the
cohomology $H^*(P_{\Sigma^{n-1}B\Z/p}X)$ is finitely generated as
an algebra over ${\mathcal A}_p$. Moreover,
$\Omega^{n-1}P_{\Sigma^{n-1}B\Z/p}X$ is weakly equivalent to
$P_{B\Z/p}\Omega^{n-1}X$, which is $B\Z/p$-local, so the induction
hypothesis applies. The cohomology of $P_{B\Z/p}X \simeq P_{B\Z/p}
P_{\Sigma^{n-1}B\Z/p}X$ is finitely generated as an algebra over
the Steenrod algebra.

Finally, since $H^*(P_{B\Z/p}X)$ is locally finite, this implies
that the space $P_{B\Z/p}X$ has finite mod $p$ cohomology.
\end{proof}

Combining this last result with Theorem~\ref{Teorema de
estructura1} we obtain:

\begin{theorem}
\label{estructura-finitamente-generados-sobre-A}
Let $X$ be a connected $H$-space such that $H^*(X)$ is a finitely
generated algebra over the Steenrod algebra. Then, $X$ is the total
space of an $H$-fibration
$$
F \rTo X \rTo Y
$$
where $Y$ is an $H$-space with finite mod $p$ cohomology and $F$
is a $p$-torsion $H$-Postnikov piece whose homotopy groups are
finite direct sums of copies of cyclic groups $\Z/p^r$ and
Pr\"ufer groups $\Z_{p^\infty}$. \hfill{\qed}
\end{theorem}

The analogous result for loop spaces follows from
Corollary~\ref{Teorema de estructura-loop spaces}.

We propose finally an extension of Hubbuck's Torus Theorem on
homotopy commutative $H$-spaces. At the prime $2$, we have:

\begin{corollary}
\label{Corolario1}
Let $X$ be a connected homotopy commutative $H$\!-space such that
the mod $2$ cohomology $H^*(X)$ is finitely generated as algebra
over the Steenrod algebra $\A_2$. Then, up to $2$-completion, $X$
is homotopy equivalent to $(S^1)^n\times F$, where $F$ is a
connected $2$-torsion $H$-Postnikov piece.
\end{corollary}

\begin{proof}
Consider the fibration $F \rightarrow X \rightarrow P_{B\Z/2} X$.
We know from the preceding theorem that the fiber is a $2$-torsion
$H$-Postnikov piece and the basis is an $H$-space with finite mod
$2$ cohomology. Both are homotopy commutative. In particular, the
mod $2$ Torus Theorem of Hubbuck, \cite{MR38:6592}, implies that
$P_{B\Z/2} X$ is, up to $2$-completion, a finite product of
circles $(S^1)^n$. Since the fiber is $2$-torsion, the above
fibration splits (not necessarily by an $H$-map) and the result
follows.
\end{proof}

When $X$ is a mod $2$ finite $H$-space, this corollary is the
original Torus Theorem due to Hubbuck. When $X$ is an $H$-space
with noetherian cohomology, $QH^*(X)\in \mathcal U_0$, the
Postnikov piece $F$ is an Eilenberg-Mac Lane space $K(P,1)$ where
$P$ is a $2$-torsion abelian group, and we get back Slack's
results \cite{MR92k:55015}, as well as their generalization by Lin
and Williams in \cite{MR92b:55011}: up to $2$-completion, $X$ is
the product of a finite number of $S^1$'s, $K(\Z/2^r,1)$'s, and
$K(\Z,2)$'s. Of course, in our setting it is no longer true that
the fiber $F$ in
Theorem~\ref{estructura-finitamente-generados-sobre-A} is a
product of Eilenberg-Mac Lane spaces.

At odd primes, there are many more finite $H$-spaces which are
homotopy commutative (all odd dimensional spheres for example).
However, Hubbuck's result still holds for finite loop spaces of
$H$-spaces, as was shown in~\cite{Aguade} by Aguad\'e and Smith.
Therefore, replacing the original Torus Theorem by the
Aguad\'e-Smith version, the same conclusion as in
Corollary~\ref{Corolario1} holds at odd primes for the loop space
on an $H$-space.


\providecommand{\bysame}{\leavevmode\hbox
to3em{\hrulefill}\thinspace}
\providecommand{\MR}{\relax\ifhmode\unskip\space\fi MR }
\providecommand{\MRhref}[2]{%
  \href{http://www.ams.org/mathscinet-getitem?mr=#1}{#2}
} \providecommand{\href}[2]{#2}



\bigskip
{\small
\begin{center}
\begin{minipage}[t]{8 cm}
Nat\`{a}lia Castellana and Jer\^{o}me Scherer\\ Departament de Matem\`atiques,\\
Universitat Aut\`onoma de Barcelona,\\ E-08193 Bellaterra, Spain
\\\textit{E-mail:}\texttt{\,natalia@mat.uab.es}, \\
\phantom{\textit{E-mail:}}\texttt{\,\,jscherer@mat.uab.es}
\end{minipage}
\begin{minipage}[t]{7 cm}
Juan A. Crespo \\ Departamento
de Econom\'\i a,\\
Universidad Carlos III de Madrid,\\
E--28903 Getafe, Spain \\
\\\textit{E-mail:}\texttt{\,jacrespo@eco.uc3m.es},
\end{minipage}
\end{center}}

\end{document}